\documentclass[12pt]{amsart}
\usepackage{amssymb}

\usepackage{amsthm}
\usepackage{amsmath}

\theoremstyle{plain}
\newtheorem{theorem}{Theorem}

\newtheorem{corollary}[theorem]{Corollary}

\newtheorem{lemma}[theorem]{Lemma}

\theoremstyle{definition}
\newtheorem{definition}[theorem]{Definition}

\numberwithin{theorem}{section}

\begin{document}
\title{On the cogeneration of cotorsion pairs}
\author{Paul C. Eklof}
\address{Math Dept, UCI\\
Irvine, CA 92697-3875}
\author{Saharon Shelah}
\address{Institute of Mathematics, Hebrew University\\
Jerusalem 91904, Israel}
\author{Jan Trlifaj}
\address{Katedra algebry MFF UK, Sokolovsk\'{a} 83, 186 75 Prague 8\\
Czech Republic}
\thanks{First author partially supported by NSF DMS-0101155. Second author supported
by the German-Israel Foundation for Scientific Research \& Development.
Publication No. 814. Third author supported by grants GA\v CR 201/03/0937
and MSM 113200007.}
\date{\today}

\begin{abstract}
Let $R$ be a Dedekind domain. In \cite{ET2}, Enochs' solution of the Flat
Cover Conjecture was extended as follows: ($*$) If $\mathfrak C$ is a
cotorsion pair generated by a class of cotorsion modules, then $\mathfrak C$
is cogenerated by a set. We show that ($*$) is the best result provable in
ZFC in case $R$ has a countable spectrum: the Uniformization Principle UP$%
^{+}$ implies that $\mathfrak C$ is not cogenerated by a set whenever $%
\mathfrak C$ is a cotorsion pair generated by a set which contains a
non-cotorsion module.
\end{abstract}

\maketitle

\section{Introduction}

For any ring $R$, if $\mathcal{S}$ is a class of (right) $R$- modules, we define 
\begin{equation*}
^{\perp }\mathcal{S}=\{A:\hbox{Ext}_{R}^{1}(A,M)=0\hbox{ for all }M\in 
\mathcal{S}\}
\end{equation*}
and 
\begin{equation*}
\mathcal{S}^{\perp }=\{A:\hbox{Ext}_{R}^{1}(M,A)=0\text{ for all }M\in 
\mathcal{S}\}
\end{equation*}
If $\mathcal{S}$ is a set (not a proper class), then $^{\perp }\mathcal{S}=\
^{\perp }\{K\}$ where $K$ is the direct product of the elements of $\mathcal{%
S}$, and $\mathcal{S}^{\perp }=\{B\}^{\perp }$ where $B$ is the direct sum
of the elements of $\mathcal{S}$ . (Henceforth, in an abuse of notation, we
will write $^{\perp }K$ instead of $^{\perp }\{K\}$, and $B^{\perp }$
instead of $\{B\}^{\perp }$.)

A \textit{cotorsion pair} (originally called a \textit{cotorsion theory}) is
a pair $\mathfrak{C}=(\mathcal{F},\mathcal{C})$ such that $\mathcal{F}={}^{\perp
}\mathcal{C}$ and $\mathcal{C}=\mathcal{F}^{\perp }$. $\mathfrak{C}$ is said to
be \textit{generated} (resp., \textit{cogenerated}) by $\mathcal{S}$ when $%
\mathcal{F}={}^{\perp }\mathcal{S}$ (resp., $\mathcal{C}=\mathcal{S}^{\perp }
$). 

A motivating example (for $R$ a Dedekind domain) is the pair $(\mathcal{F},%
\mathcal{C})$ where $\mathcal{F}$ is the class of torsion-free modules and $%
\mathcal{C}=\mathcal{F}^{\perp }$; the members of $\mathcal{C}$ are called
cotorsion modules. Equivalently, $K$ is cotorsion if and only if $\hbox{%
Ext}_{R}^{1}(Q,K)=0$, where $Q$ is the quotient field of $R$ (cf. \cite[\S XIII.8]{FS}. Pure-injective
modules are cotorsion, and torsion-free cotorsion modules are pure-injective.

Cotorsion theories were first studied by Salce \cite{Sal}; their study was
given new impetus by the work of G\"{o}bel-Shelah \cite{GS}. (See, for
example, \cite[Chap. XVI]{EM} for an introduction to these concepts.)

\smallskip\ 

In this paper we are interested in the question of when a cotorsion pair $(%
\mathcal{F},\mathcal{C})$ is cogenerated by a set, or, equivalently, when
there is a single module $B\in \mathcal{F}$ such that $\mathcal{C}=B^{\perp
} $. One reason this question is of interest is that, by a result in \cite
{ET1}, if $(\mathcal{F},\mathcal{C})$ is cogenerated by a set, then it is 
\textit{complete}, that is, for every module $M$, there is an epimorphsim $\psi: N \to M$ such that $N \in \mathcal{F}$ and ker$(\psi) \in \mathcal{C}$; in
particular, $\mathcal{F}$-precovers exist for all $R$-modules. It is these
ideas and results that are involved in the proof of the Flat Cover
Conjecture by Enochs \cite{BBE}; see the introduction to \cite{ET2} for the
historical sequence of events. (See also \cite{EJ} and/or \cite{Xu} for a
comprehensive study of (pre)covers and their uses.)

The following is proved in \cite{ET2}:

\begin{theorem}
\label{covers}For any ring $R$, if $\mathfrak{C}=(\mathcal{F},\mathcal{C})$ is a
cotorsion pair which is generated by a class of pure-injective modules, then 
$\mathfrak{C}$ is cogenerated by a set. Moreover, if $R$ is a Dedekind domain,
the same conclusion holds when $\mathfrak{C}$ is generated by a class of
cotorsion modules, or, equivalently, when every element of $\mathcal{C}$ is
cotorsion.
\end{theorem}

Note that $(\mathcal{F},\mathcal{C)}\ $is generated by a class of cotorsion
modules if and only if $Q\in \mathcal{F}$, in which case every member of $%
\mathcal{C}$ is cotorsion.

\smallskip\ 

The case when $\mathcal{C}$ contains non-cotorsion modules is more
complicated, and the results depend on the extension of ZFC we work in. In 
\cite{ET2} it is proved that it is consistent with ZFC that the conclusion
of Theorem \ref{main} holds for even more cotorsion pairs:

\begin{theorem}
\label{hered}G\"{o}del's Axiom of Constructibility (V = L) implies that $%
\mathfrak C$ is cogenerated by a set whenever $\mathfrak C$ is a cotorsion
pair generated by a set and $R$ is a right hereditary ring.
\end{theorem}

\smallskip\ 

The main result of this paper is that Theorem \ref{covers} is the best that
can be proved in ZFC (even in ZFC + GCH) for cotorsion pairs which are
generated by a set --- at least for certain rings, including $\mathbb{Z}$:

\begin{theorem}
\label{main} It is consistent with ZFC + GCH that if $R$ is a Dedekind
domain with a countable spectrum and $\mathfrak C=(\mathcal{F},\mathcal{C})$
is a cotorsion pair generated by a set which contains a non-cotorsion
module, then $\mathfrak C$ is \emph{not} cogenerated by a set.
\end{theorem}

The assumption that $\mathfrak C$ is generated by a set is essential in \ref
{main}: for example, by a classical result of Kaplansky, the cotorsion pair $%
(\mathcal{P}_{0},\mbox{Mod-}R)$ is cogenerated by a set (of countably
generated modules), for any ring $R$. (Here, $\mathcal{P}_{0}$ denotes the
class of all projective modules.)

Putting together Theorems \ref{covers} and \ref{main}, we have:

\begin{corollary}
\label{maincor} Let $R$ be a Dedekind domain with a countable spectrum, and
let $K$ be an $R$-module. It is provable in ZFC + GCH that there is a module 
$B$ such that $(^{\perp }K)^{\perp }=B^{\perp }$ if and only if $K$ is
cotorsion.
\end{corollary}

\noindent \textsc{Proof}. If $K$ is cotorsion, it is proved in \cite{ET2}
that $B$ exists. (This is provable in ZFC alone.) The other direction
follows immediately from Theorem \ref{main} for the cotorsion pair $(^{\perp
}K,(^{\perp }K)^{\perp })$. \qed  

\smallskip\ 

In \cite{ES2} this result was proved for \textit{countable torsion-free }$%
\mathbb{Z}$-modules $K$. It was also proved there that the cotorsion pair $%
(^{\perp }\mathbb{Z},(^{\perp }\mathbb{Z})^{\perp })$ is not complete.

\smallskip\

Theorem \ref{main} is proved in the next two sections. In the first one we prove in ZFC some preliminary results. In the following section we invoke the additional set-theoretic hypothesis UP$^+$.

\section{Results in ZFC}

We will make use of the following result from \cite{ET1}. (See also 
\cite[XVI.1.2 and XVI.1.3]{EM}.)

\begin{theorem}
\label{vanext}Let $B$ be an $R$-module and let $\kappa $ be a cardinal 
 $ > |R|+|B|$. Let $\mu $ be a cardinal $>\kappa $ such that $\mu
^{\kappa }=\kappa $ . Then there is a module $A\in B^{\perp }$ such that $%
A=\bigcup_{\nu <\mu }A_{\nu }$ (continuous), $A_{0}=0$ (or any given module
of size $<\kappa $), and such that for all $\nu <\mu $, $A_{\nu +1}/A_{\nu }$
is isomorphic to $B$.

Moreover, if, for some $R$-module $K$, $B\in {}^{\perp }K$, then  $%
A/A_{\nu } \in {}^{\perp }K$ for all $\nu <\mu $.\qed  
\end{theorem}

The continuity condition on the $A_{\nu }$ means that for every limit
ordinal $\sigma <\mu $, $A_{\sigma }=\bigcup_{\nu <\sigma }A_{\nu }$.

\smallskip

\relax From now on, $R$ will denote a Dedekind domain and $Q$ will denote its
quotient field. Moreover, we assume that $Q$ is countably generated as an $R$%
-module, or, equivalently, that $R$ has a countable spectrum.

The conditions on $A$ in Theorem \ref{vanext} motivate the hypotheses in the
following lemmas. Recall that a module $M$ is \textit{reduced} if $\mbox{Hom}%
_{R}(Q,M)=0$.

\begin{lemma}
\label{noQ}Let $B$ be a torsion-free reduced module. Let $\mu $ be a limit
ordinal and suppose $M=\bigcup_{\nu <\mu }M_{\nu }$ (continuous), where $%
M_{0}=0$, and for all $\nu <\mu $, $M_{\nu +1}/M_{\nu }$ is isomorphic to $B$%
. Then $M$ is torsion-free and reduced.
\end{lemma}

\noindent \textsc{Proof}. It is clear that $M$ is torsion-free. Suppose that
there is a non-zero homomorphism, hence an embedding, $\theta :Q\rightarrow
M $. Let $\tau $ be minimal such that $M_{\tau }$ contains a non-zero
element, $\theta (y)$, of the range of $\theta $. Then $\tau $ is not a
limit ordinal; say $\tau =\nu +1$, and $\theta $ induces a non-zero map,
hence an embedding, of $Q$ into $M/M_{\nu }$. Since $M/M_{\nu +1}$ has no
torsion, this map embeds $Q$ into $M_{\nu +1}/M_{\nu }$, which is a
contradiction, since $M_{\nu +1}/M_{\nu }\cong B$. \qed  

\smallskip\

\begin{definition}
\label{rho} By hypothesis on $R$ we can fix a countable set $\{\rho
_{j}:j\in \omega \}$ of non-units of $R$ such that $\{(\prod_{i<j}\rho
_{i})^{-1}:j\in \omega \}$ generates $Q$ as an $R$-module.
\end{definition}

\begin{lemma}
\label{nosol}Let $B$ be a torsion-free $R$-module. Suppose $M=\bigcup_{n\in
\omega }M_{n}$ such that $M_{0}=0$, and for all $n\in \omega $, $%
M_{n+1}/M_{n}$ is isomorphic to $B$. Suppose that for some $k\in \omega $
and all $n\in \omega $, $a_{n}+M_{n}$ is an element of $M_{n+1}/M_{n}$ which
does not belong to $\rho _{k}(M_{n+1}/M_{n}).$ Then the system of equations 
\[
\{\rho _{n}v_{n+1}=v_{n}-a_{n}:n\in \omega \}
\]
in the variables $\{v_{n}:n\in \omega \}$does not have a solution in $M$.
\end{lemma}

\noindent \textsc{Proof}. Suppose, to the contrary, that there is a solution 
$v_{n}=u_{n}\in M$. We have $u_{0}\in M_{m}$ for some $m\geq k$. Since $%
a_{n}\in M_{m}$ for $n<m$, and since $B$ is torsion-free, $u_{n}\in M_{m}$
for $n\leq m$. But then $\rho _{m}u_{m+1}=u_{m}-a_{m}$ implies that $%
u_{m+1}+M_{m}$ belongs to $M_{m+1}/M_{m}$ (since $M/M_{m+1}$ is
torsion-free) and thus $\rho _{k}$ divides $a_{m}+M_{m}$ in $M_{m+1}/M_{m}$,
which contradicts the choice of $a_{m}$. \qed  

\smallskip\ 

Recall that a module $M$ is called a \textit{splitter }if $\hbox{Ext}%
_{R}^{1}(M,M)=0$. (See, for example, \cite{Sch}, \cite{GS}, or \cite[Chap. XVI]{EM}.)

\begin{lemma}
\label{tfsplitter} If $\mathfrak C$ is a cotorsion pair which is generated
and cogenerated by sets, then there is a torsion-free splitter which
generates $\mathfrak C$.
\end{lemma}

\noindent \textsc{Proof}. Let $\mathfrak C=(\mathcal{F},\mathcal{C})$. Let $%
B,K$ be modules such that $\mathcal{F}={}^{\perp }K$ and $\mathcal{C}%
=B^{\perp }$. By \cite[Theorem 10]{ET1}, $K$ has a special $\mathcal{F}$%
-precover, i.e., there is an exact sequence $0\to M\to N\to K\to 0$ such
that $M\in \mathcal{C}$ and $N\in \mathcal{F}$. Since $K\in \mathcal{C}$,
also $N\in \mathcal{C}$, and $N\in \mathcal{C}\cap \mathcal{F}$ is a
splitter.

We have $\mathcal{F}={}^{\perp }N$ (since clearly $\mathcal{F}\subseteq
{}^{\perp }N$, and $^{\perp }N\subseteq {}^{\perp }K=\mathcal{F}$). Let $T$
be the torsion part of $N$. Then $T$ is a direct sum of its $p$-components, $%
T=\bigoplus_{p\in mSpec(R)}T_{p}$ . If $T_{p}\neq 0$, then $\mbox{Ext}%
_{R}^{1}(R/p,N)=0$, so $\mbox{Hom}_{R}(R/p,E(N)/N)=0$, and hence $\mbox{Hom}%
_{R}(R/p,E(T_{p})/T_{p})=0$. Therefore $T_{p}$ is divisible. So $N=T\oplus L$
where $L$ is a torsion-free splitter. Since $T$ is divisible, $^{\perp
}L={}^{\perp }N=\mathcal{F}$. \qed

\smallskip\

\begin{lemma}
\label{Harrison} Suppose that $\mathfrak C$ is a cotorsion pair which is
cogenerated by a cotorsion module, and generated by a set. Then $\mathfrak C$
is cogenerated by a cotorsion module of the form $B\oplus T$ where $B$ is
torsion-free, $T$ is torsion, and for every prime $p$ such that $R/p$ is a
submodule of $T$, $pB=B$.
\end{lemma}

\noindent \textsc{Proof}. Let $\mathfrak C=(\mathcal{F},\mathcal{C})$ and
let $K$ be a module such that $\mathcal{F}={}^{\perp }K$. If $K$ is
cotorsion, then by \cite[Thm. 16]{ET2}, there is a set of maximal ideals $P$
such that $\mathcal{F}$ is the set of all modules with zero $p$-torsion part
for all $p\in P$. Then $\mathcal{C}=B^{\perp }$ where $B=Q\oplus
\bigoplus_{q\notin P}R/q$.

So we can assume that $K$ is not cotorsion, and that, by Lemma \ref
{tfsplitter}, $K$ is torsion-free.

Let $C$ be a cotorsion module such that $\mathcal{C}=C^{\perp }$. We have $%
C=D\oplus E$ where $D$ is divisible and $E$ reduced. Since $K$ is not
cotorsion, $D$ is torsion. Denote by $T^{\prime }$ the torsion part of $E$.
By a theorem of Harrison-Warfield, \cite[XIII.8.8]{FS}, we have $E=B\oplus G$
where $B$ is torsion-free reduced and pure-injective, and $G$ is a cotorsion
hull of $T^{\prime }$. We claim that there is an exact sequence $0\to
T^{\prime }\to G\to Q^{(\delta )}\to 0$ for some $\delta \geq 0$. 

Indeed,  by \cite[3.4.5]{Xu}, $G$ is a cotorsion envelope of $T^{\prime }$
in the sense of Enochs. Now by Theorem \ref{vanext} there is a cotorsion
preenvelope $G^{\prime }$ of $T^{\prime }$ such that $G^{\prime }/T^{\prime }
$ is the union of a continuous chain with successive quotients 
isomorphic to $Q$, and hence $G^{\prime }/T^{\prime }\cong Q^{(\gamma )}$
for some $\gamma $. The claim now follows since $G/T^\prime$ is isomorphic 
to a direct summand of $G^{\prime }/T^\prime$ by \cite[1.2.2]{Xu}

Since $K$ is torsion-free and $G\in \mathcal{C}$, an application of $%
\mbox{Hom}_{R}(-,K)$ yields

\[
0=\mbox{Hom}_{R}(T^{\prime },K)\to \mbox{Ext}_{R}^{1}(Q^{(\delta )},K)\to %
\mbox{Ext}_{R}^{1}(G,K)=0. 
\]
Thus, $\hbox{Ext}(Q^{(\delta )},K)=0$, so since $K$ is not cotorsion, $%
\delta =0$ and $T^{\prime }=G$. Hence $C=B\oplus T$ where $T=T^{\prime
}\oplus D$ is torsion.

By \cite[5.3.28]{EJ}, there is a set $P$ of maximal ideals of $R$ such that $%
B\cong \prod_{p\in P}J_{p}$ where $J_{p}$ is the $p$-adic completion of a
free module over the localization of $R$ at $p$. In particular, $qB=B$ for
all maximal ideals $q\notin P$. For each $p\in P$, there is an exact
sequence $0\to J_{p}\to E(J_{p})\to I_{p}\to 0$ where $I_{p}$ is a direct
sum of copies of $E(R/p)$, and $E(J_{p})=Q^{(\alpha _{p})}$ for some $\alpha
_{p}>0$.

Let $q$ be a maximal ideal such that $R/q$ embeds in $T$. Assume $q\in P$.
Then an application of $\mbox{Hom}_{R}(-,K)$ yields

\[
0=\mbox{Ext}_{R}^{1}(I_{q},K)\to \mbox{Ext}_{R}^{1}(Q^{(\alpha _{q})},K)\to %
\mbox{Ext}_{R}^{1}(J_{q},K)=0\text{.}
\]
The first Ext is zero because $R/q\hookrightarrow T$;  so $R/q\in \mathcal{F} = {}^{\perp }\mathcal{C}$ and thus $E(R/q)\in \mathcal{F}$ by \cite[Lemma 1]
{ET1} since $E(R/q)$ is the union of a continuous chain of modules with
successive quotients isomorphic to $R/q$; the last Ext is zero because $%
J_{q}\in \mathcal{F}$. So $K$ is cotorsion, a contradiction. This proves
that $q\notin P$ and hence $qB=B$. \qed   

\section{Proof of Theorem \ref{main}}

Let $\mathfrak C = (\mathcal{F},\mathcal{C})$ be a cotorsion pair
cogenerated by a set, and generated by a non-cotorsion module $K$. We aim to
produce a contradiction by constructing $H\in {}^{\perp }K$ ($=\mathcal{F}$)
and $A\in \mathcal{C}$ such that $\hbox{Ext}^1_R(H,A)\neq 0$. We do this
assuming GCH plus the following principle, which is consistent with ZFC +
GCH (cf. \cite{ES} or \cite{SS}):

\begin{quote}
\textbf{(UP}$^{+}$\textbf{)} For every cardinal $\mu $ of the form $\tau
^{+} $ where $\tau $ is singular of cofinality $\omega $ there is a
stationary subset $S$ of $\mu $ consisting of limit ordinals of cofinality $%
\omega $ and a ladder system $\bar{\zeta}=\{\zeta _{\delta }:\delta \in S\}$
which has the $\lambda $-uniformization property for every $\lambda <\tau $.
\end{quote}

\noindent Recall that if $S$ is a subset of an uncountable cardinal $\mu $
which consists of ordinals of cofinality $\mathcal{\omega }$, a\textit{\
ladder system} on $S$ is a family $\bar{\zeta}=\{\zeta _{\delta }:\delta \in
S\}$ of functions $\zeta _{\delta }:\mathcal{\omega }\rightarrow \delta $
which are strictly increasing and have range cofinal in $\delta $. For a
cardinal $\lambda $, we say that $\bar{\zeta}$ has the $\lambda $\textit{%
-uniformization property} if for any functions $c_{\delta }:\mathcal{\omega }%
\rightarrow \lambda $ for $\delta \in S$, there is a pair $(f,f^{*})$ where $%
f:\mu \rightarrow \omega $ and $f^{*}:S\rightarrow \mathcal{\omega }$ such
that for all $\delta \in S$, $f(\zeta _{\delta }(\nu ))=c_{\delta }(\nu )$
whenever $f^{*}(\delta )\leq \nu <\mathcal{\omega }$. \noindent We refer to 
\cite[Chap. XIII]{EM} for more details.

\smallskip\ 

We consider two cases: (1) $\mathfrak C$ is cogenerated by a cotorsion
module; and (2) the negation of (1).

The module $H$ will be the same in both cases. Let $\bar{\zeta}=\{\zeta
_{\delta }:\delta \in S\}$ be as in (UP$^{+}$) for this $\mu $. We also use
the notation from Definition \ref{rho}. Let $H=F/L$ where $F$ is the free
module with the basis $\{y_{\delta ,n}:$ $\delta \in S$, $n\in \omega \}\cup
\{x_{j}:$ $j<\mu \}$ and $L$ is the free submodule with the basis $%
\{w_{\delta ,n}:$ $\delta \in S$, $n\in \omega \}$ where 
\begin{equation}
w_{\delta ,n}=y_{\delta ,n}-\rho _{n}y_{\delta ,n+1}+x_{\zeta _{\delta }(n)}.
\label{eqnsG}
\end{equation}

Then $H$ is a module of cardinality $\mu $ and the uniformization property
of $\bar{\zeta}$ implies that $H$ $\in {}^{\perp }K$. (In fact, $H$ $\in
{}^{\perp }K$ for any module $K$ of cardinality $<\tau $. See 
\cite[Chap. XIII]{EM} or \cite{T}.)

\smallskip\ 

Assuming we are in Case (1), let $B\oplus T$ be a cogenerator of $\mathfrak %
C $ as given in Lemma \ref{Harrison}. Let $\kappa \geq \max (|B|,$ $|R|,|K|)$
and let $\mu =\tau ^{+}=2^{\tau }$ where $\tau >\kappa $ is a singular
cardinal of cofinality $\omega $. Then $\mu ^{\kappa }=\mu $. Let $%
A=\bigcup_{\nu <\mu }A_{\nu }$ be as in Theorem \ref{vanext} for this $B$
and $\mu $; so, in particular, $A\in B^{\perp }$. Note that then $A\in
(B\oplus T)^{\perp }=\mathcal{C}$ because $T^{\perp }$ consists of precisely
those modules $M$ such that $pM=M$ whenever $R/p\hookrightarrow T$. Note
that $A/A_{\delta }$ is torsion-free for all $\delta \in \mu $, because $B$
is torsion-free.

We need to show that $\hbox{Ext}^1_R(H,A)\neq 0$; in other words, to
define a homomorphism $\psi :L\rightarrow A$ which does not extend to $F$.

Since $B$ is reduced there is a $k\in \omega $ such that $\rho _{k}B\neq B$;
then for all $\delta \in S$ and $n\in \omega $ we can choose $a_{\delta
,n}\in A_{\delta +n+1}$ such that $a_{\delta ,n}+A_{\delta +n}\notin \rho
_{k}(A_{\delta +n+1}/A_{\delta +n})$. We claim that

\begin{quotation}
($\maltese $) for all $\delta \in S$, the family of equations 
\[
\mathcal{E}_{\delta }=\{\rho _{n}v_{n+1}=v_{n}-(a_{\delta ,n}+A_{\delta
}):n\in \omega \} 
\]
does not have a solution in $A/A_{\delta }$.
\end{quotation}

Supposing, for the moment, that this claim is true, we will prove that $%
\hbox{Ext}^1_R(F/L,A)\neq 0$. Define $\psi :L\rightarrow A$ by $\psi
(w_{\delta ,n})=a_{\delta ,n}$ for all $\delta \in S$, $n\in \omega $.
Suppose, to obtain a contradiction, that $\psi $ extends to a homomorphism $%
\varphi :F\rightarrow A$. The set of $\delta <\mu $ such that $\varphi
(x_{j})\in A_{\delta }$ for all $j<\delta $ is a club, $C$, in $\mu $, so
there exists $\delta \in S\cap C$. By applying $\varphi $ to the relations (%
\ref{eqnsG}), and since $\varphi (x_{j})\in A_{\delta }$ for all $j<\delta $%
, we have that $v_{n}=\varphi (y_{\delta ,n})+A_{\delta }$ is a solution to
the equations in $A/A_{\delta }$, a contradiction.

Thus it remains to prove ($\maltese $). Suppose that ($\maltese $) is false
for some $\delta \in S$, and that for some $\{b_{n}:n\in \omega \}\subseteq
A $, $v_{n}=b_{n}+A_{\delta }$ is a solution to $\mathcal{E}_{\delta }$.
There are two subcases.

Suppose first that $b_{0}+A_{\delta +\omega }$ is a non-zero element of $%
A/A_{\delta +\omega }$. Then $A/A_{\delta +\omega }$ contains a copy of $Q$
(generated over $R$ by the cosets of the $b_{n}$, $n\in \omega $). But this
contradicts Lemma \ref{noQ} (with $M=A/A_{\delta +\omega }$, $M_{\nu
}=A_{\delta +\omega +\nu }/A_{\delta +\omega }$).

Otherwise we can prove by induction that $b_{n}\in A_{\delta +\omega }$ for
all $n\in \omega $ because $A/A_{\delta +\omega }$ has no torsion and $\rho
_{n}(b_{n+1}+A_{\delta +\omega })=b_{n}+A_{\delta +\omega }$. Thus there is
a solution of 
\[
\{\rho _{n}v_{n+1}=v_{n}-(a_{\delta ,n}+A_{\delta }):n\in \omega \} 
\]
in $A_{\delta +\omega }/A_{\delta }$. But this contradicts Lemma \ref{nosol}
(with $M=A_{\delta +\omega }/A_{\delta }$, $M_{n}=A_{\delta +n}/A_{\delta }$
and $a_{n}=a_{\delta ,n}+A_{\delta }$).

This completes the proof in Case (1).

\smallskip\ 

Now supposing we are in Case (2), let $B$ be a module cogenerating $%
\mathfrak C$. Let $\kappa \geq \max (|B|,$ $|R|,|K|)$ and let $\mu =\tau
^{+}=2^{\tau }$ where $\tau >\kappa $ is a singular cardinal of cofinality $%
\omega $. Let $A=\bigcup_{\nu <\mu }A_{\nu }$ be as in Theorem \ref{vanext}
for this $B$ and $\mu $; so $A\in B^{\perp }$. Let $H$ be as above.

Then for all $\delta \in \mu $, $A/A_{\delta }$ cogenerates $\mathfrak C$
since the construction of $B$ and Lemma 1 of \cite{ET1} implies that $M\in
(A/A_{\delta })^{\perp }$ whenever $M\in B^{\perp }$. Hence, since we are in
Case (2), $\hbox{Ext}_{R}^{1}(Q,A/A_{\delta })\neq 0$ for all $\delta \in
\mu $.

Now $Q\cong F_{\delta }/L_{\delta }$ where $F_{\delta }$ is the free module
with the basis $\{y_{\delta ,n}:$ $n\in \omega \}$ and $L_\delta$ is the
free submodule with the basis $\{w_{\delta ,n}^{\prime }:$ $\delta \in S$, $%
n\in \omega \}$ where $w_{\delta ,n}^{\prime }=y_{\delta ,n}-\rho
_{n}y_{\delta ,n+1}$. Hence there is a homomorphism $\psi _{\delta
}:L_{\delta }\rightarrow A/A_{\delta }$ which does not extend to $F_{\delta
} $.

Let $\pi _{\delta }:A\to A/A_{\delta }$ be the canonical projection. Define $%
\psi :L\rightarrow A$ so that $\pi _{\delta }\psi (w_{\delta ,n})=\psi
_{\delta }(w_{\delta ,n}^{\prime })$. In order to prove $\mbox{Ext}%
_{R}^{1}(H,A)\neq 0$, we will show that $\psi $ does not extend to a
homomorphism $\varphi :F\rightarrow A$. If it did, there would exist $\delta
\in S\cap C$ where $C$ is the club of all $\delta <\mu $ such that $\varphi
(x_{j})\in A_{\delta }$ for all $j<\delta $. But then $\pi _{\delta }\circ
(\varphi \upharpoonright F_{\delta })$ would be an extension of $\psi
_{\delta }$, a contradiction.

This completes the proof of Theorem \ref{main}. \qed

\end{document}